\documentclass[12pt,a4paper]{article}
\usepackage{amsmath, ams symb, theorem, latexsym}
\newtheorem{theorem}{Theorem}[section] 

\newtheorem{proposition}[theorem]{Proposition}

\newtheorem{corollary}[theorem]{Corollary}

\newtheorem{remark}[theorem]{Remark}

\allowdisplaybreaks

\newcommand\finbox{~\hfill$\Box$}%
\begin{document}
\def\Om {{\Omega}}
\def\la {{\lambda}}
\def\ep {{\epsilon}}
\def \zz { {\mathbb Z}}
\def \qz { {\mathbb Q}}
\def \rz { {\mathbb R}}
\def\Ag {{\cal A}} 
\def\Bg {{\cal B}} 
\def\Cg {{\cal C}} 
\def\Zg {{\cal Z}} 
\def\Fg {{\cal F}} 
\def\Hg {{\cal H}} 
\def\Hb {{\bf H}} 
\def\hb {{\bf h}} 
\def\Ig {{\cal I}} 
\def\Jg {{\cal J}} 
\def\Eg {{\cal E}} 
\def\Kg {{\cal K}} 
\def\Lg {{\cal L}} 
\def\Ng {{\cal N}} 
\def\Og {{\cal O}} 
\def\Pg {{\cal P}} 
\def\Qg {{\cal Q}} 
\def\Rg {{\cal R}} 
\def\Rb {{\bf R}} 
\def\Nb {{\bf N}} 
\def\Sg {{\cal S}} 
\def\Vg {{\cal V}} 
\def\Wg {{\cal W}} 
\def\Ug {{\cal U}} 
\def\qz {{\mathbb Q}}
\def \zz {{\mathbb Z}}
\def \rz {{\mathbb R}}
\def \cz {{\mathbb C}}
\def \dist {{\rm dist\,}} 
\def \Inte{{\rm Int\,}}
\def \Cap{{\rm Cap\,}} 
\def \vph{{\varphi}}
\newcommand {\ar}{\rightarrow}
\newcommand {\pa}{\partial}
\newcommand {\oa}{\overrightarrow}
\numberwithin{equation}{section} 
\centerline
{{\bf{\large  Eigenfunctions for rectangles with Neumann boundary conditions}}}

\centerline{}

\centerline{}

\centerline{Thomas Hoffmann-Ostenhof, University of Vienna}
\centerline{}
{{
\centerline{\today}
\begin{abstract}

 Consider the eigenfunctions $u$ for a free rectangular membrane so that $-\Delta u=\la u$
on $\mathcal R(c,d)=(0,c)\times (0,d)$. In this note we show that if  $u>0$  on $\pa \mathcal R(c,d)$ then
$u\equiv C>0$ for some positive constant $C$.       
\end{abstract} 

\section{ Introduction and the result}  
We consider the Laplacian on  a rectangle $\mathcal R(c,d)=(0,c)\times (0,d)$ with Neumann boundary conditions so that 
$-\Delta u=\la u$ and $\frac{\pa}{\pa\nu}u=0$ on $\pa \mathcal R(c,d)$. 
The Neumann eigenvalues are given by 
\begin{equation}\label{lamn} 
\la_{m,n}=\pi^2(m^2/c^2+n^2/d^2),\: (m,n)\in \mathbb N_0^2
\end{equation}
and the associated eigenfunctions are given, up to a multiplicative constant, by 
\begin{equation}\label{umnirr}
u_{m,n}(x,y)=\cos (m\pi x/c)\cos (n\pi y/d).
\end{equation}
Of course if an eigenvalue $\la_{m,n}=\la_{m',n'}$  is not simple 
then we have to consider linear combinations of these eigenfuctions in its eigenspace.

\begin{theorem}\label{I}
Assume that $u$ is an eigenfunction on $\mathcal R(c,d)$ so that $-\Delta u=\la u$ and that 
$u$ satisfies Neumann boundary conditions, i.e. $\frac{\pa}{\pa \nu} u=0$ on $\pa \mathcal R(c,d)$ where $\nu$ denotes the 
outward directed normal.\\  
If
\begin{equation}\label{statement}
u>0 \text{ on } \pa \mathcal R(c,d) \text{ then } u\equiv \text{const}. 
\end{equation}
\end{theorem}
\begin{remark}\label{1}
Here are some immediate observations.\\  
{\bf(i)} Note that there are eigenfunctions, e.g. for the square  membrane which are positive on the boundary, 
but not strictly positive. Take for instance on $\mathcal R(2\pi ,2\pi)$ the function   
\begin{equation*}
 u=\cos x+\cos y 
\end{equation*} 
then on the boundary $u(x,0)=u(x,2\pi)=\cos x+1$ and similarily for $u(0,y)=u(2\pi,y)$.\\
For the Dirichlet case the situation is different. There is an example in Courant Hilbert \cite{CH1}, page 302,
for an eigenfunction associated to an excited eigenvalue for the square   whose zeroset does not hit the boundary.\\  
{\bf(ii)}
In corollary \ref{TorCy} we give two other examples where theorem \ref{I} holds; for the disk it
does not hold and probably it holds 
only for very special cases.  So it  would be interesting to find other domains for which theorem \ref{I} holds.\\  
{\bf (iii)} In Oberwolfach during the workshop {\em Geometric Aspects of Spectral Theory} in July 2012 I have raised the question 
whether theorem \ref{I} holds, \cite{OW}.\\
{\bf (iv)} I was not able to find results in the 
spirit of the present work. For eigenfunctions associated to high eigenvalues 
whose nodal lines hit the boundary there is recent  work, \cite{TothZelditch}. 
\\
{\bf (v)} Theorem \ref{I} is the general result, but one can asked more detailed questions. If one requires only 
in Theorem \ref{I} $u\ge 0$ on $\pa \mathcal R(c,d)$ then there are cases for which this also implies $u\equiv C>0$.
We shall not discuss this further here. Also the question whether Theorem \ref{I} extends in to cuboids is natural, 
but will not be addressed here.   
\end{remark}
{\bf Acknowledgement} I want to thank Bernard Helffer for a careful reading of the manuscript  and helpful
remarks and Nikolai Nadirashvili for motivating discussions.  
 
\section{Proof}
We assume for contradiction that there is an eigenfunction which is strictly positive on the boundary of our rectangle. 
We have to consider three cases.   

\begin{equation}\label{C1}
\frac{c^2}{d^2}\not\in \mathbb Q, 
\end{equation}
\begin{equation}\label{C2} 
\frac{c}{d}\not\in \mathbb Q,\text{ but }\frac{c^2}{d^2}\in\mathbb Q 
\end{equation} 
\begin{equation}\label{C3}
\frac{c}{d}\in \mathbb Q
\end{equation}
We start with case \eqref{C1}. \\  
   
If $c^2/d^2\not\in \mathbb Q$ then the eigenvalues
 are simple as can be seen form \eqref{lamn}. The  real 
valued eigenfunctions are given up to a multiplicative constant by 
\begin{equation}\label{umnirr}
u_{m,n}(x,y)=\cos (m\pi x/c)\cos (n\pi y/d).
\end{equation} 
Hence $u(x,0)=\cos(m\pi x/c)$ and $u(0,y)=\cos(n\pi y/d)$ and if $m>0$, respectively $n>0$ then $u$ 
has both signs on $\pa\mathcal R(c,d)$. But if $m=n=0$\ then $u$ is a constant. This settles the first case.
 \\\\
We continue with the second case \eqref{C2}. Now we  can have degenerate eigenvalues.  
Pick a non-constant eigenfunction $u$ associated to a degenerate eigenvalue $\la$ and 
assume that $u(x,0),u(0,y), u(x,d),u(c,y)$ are strictly positive. 
Consider the eigenspace $U(\la)$ and denote by 
\begin{equation}\label{Ila}
 I(\la)=\big\{(m,n)\in \mathbb N_0^2\:\big|\:\pi^2(\frac{m^2}{c^2}+\frac{n^2}{d^2})=\la\big\}.
\end{equation}
 The  eigenfunctions $u\in U(\la)$ are then given by
\begin{equation}\label{uila}
 u=\sum_{(m,n)\in I(\la)}a_{m,n}\cos(m\pi x/c)\cos(n\pi y/d) 
\end{equation}
for $a_{i,j}\in \mathbb R$. 
Consider 
\begin{equation}\label{ux0} 
 u(x,0)=\sum_{(m,n)\in I(\la)}a_{m,n}\cos(m\pi x/c)\\
\end{equation} 
and
\begin{equation}\label{u0y}
u(0,y)=\sum_{(m,n)\in I(\la)}a_{m,n}\cos(n\pi y/d).
\end{equation}
{\bf Assertion:}\\
 If $u(x,y)>0$ on $\pa \mathcal R(c,d)$  then there is a pair $(m_0,n_0)\in \mathbb N^2$ such that $(m_0,0)\in I(\la)$
 and an $(0,n_0)\in I(\la)$ so that 
\begin{equation}\label{0m0n}
m_0^2/c^2=n_0^2/d^2=\la.
\end{equation}
{\bf Proof.}  Suppose for contradiction there is no such pair. Then $u(x,0)$ is orthogonal to any constant 
so that $\int_0^cu(x,0)dx=0$ by the orthogonality of the cosines and must have both signs.  The same 
argument hold for $u(0,y)$, proving the assertion. \finbox
\\  
So \eqref{0m0n} implies  $c/d=m_0/n_0$ contradicting the assumption that 
$c/d$ is irrational. This settles the second case.    
\\\\
We continue with the third case, \eqref{C3}. 
After a suitable scaling we have a rectangle $\mathcal R(p,q)$ for relatively prime integers
$p,q$. For contradiction we assume that there is a non-constant eigenfunction $v>0$ on $\pa \mathcal R(p,q)$. 
We can symmetrize 
$v$ with respect to the two axes of symmetry $x=p/2$ and $y=q/2$  and consider $w=(v(x,y)+v(p-x,y)+v(x,q-y)+v(p-x,q-y))/4$.
Obviously also $w>0$ on $\pa \mathcal R(p,q)$. Furthermore since $\cos (m\pi x)/p$ for $x=0$ and for $x=p$ cancel for odd $m$
it suffices to consider  even $m,n$ in the expansion in cosines of $w$ in \eqref{uila} (where $c,d$ is replaced by $p,q$). 
We can now periodically continue in the $x$ and $y$ direction until we get a square with side length of the least common 
multiple of $p,q$. This square we scale to the square $Q$ with side length $2\pi$.    

The corresponding eigenfunctions are given by 
\begin{equation}\label{uQ}
u=\sum_{(m,n)\in I(\la)}a_{m,n}f_{m,n} 
\end{equation}
with $f_{m,n}=\cos mx\cos n y$. Note that also $\cos kx\cos \ell y$ for $k,\ell$ half integers would be also Neumann
eigenfunctions but those terms have been 
eliminated by the symmetrization in the preceeding paragraph.  
$I(\la)$ is defined as above as 
\begin{equation}\label{IlaQ}
I(\la)=\{(m,n)\in \mathbb N_0^2\:|\:m^2+n^2=\la\}. 
\end{equation}
Finally we can require that $u(x,y)=u(y,x)$ so that $a_{m,n}=a_{n,m}$ and $u(x,y)$ is invariant with respect to 
the symmetry operations leaving the square invariant.     
\\      
We define
$$\mathcal I=\{(i,j)\in \mathbb N^2\:|\: ij\text{ is odd }\} $$
and 
$$\mathcal J=\{(i,j)\in \mathbb N_0^2\:|\: i+j\text{  is odd } \}.$$

Obviously $(i,j)\in \mathcal I$ if and only if both $i$ and $j$ are odd, and in $\mathcal J$
if and only if one is even and the other one is odd.   
   
We continue with easy considerations:  
Clearly each positive integer $\la$ can be written either as  
\begin{equation}\label{4s}
\la =4^s(2\ell+1) 
\end{equation}
or as  
\begin{equation}\label{2.4s}
 \la=2\cdot4^s(2\ell+1)
\end{equation}
where $(s,\ell)\in \mathbb N_0^2$. 
\begin{proposition}\label{la=}
Suppose that $\la=m^2+n^2$ with $m+n>0$. If 
$\la$ satisfies \eqref{4s} then  
\begin{equation}\label{J} 
2^s\mid m,\: 2^s\mid n\text{ and }\Big(\frac{m}{2^s},\frac{n}{2^s}\Big)\in \mathcal J.
\end{equation}
If $\la$ satisfies \eqref{2.4s} then 
\begin{equation}\label{Ii}\
2^s\mid m,\: 2^s\mid n\text{ and }\Big( \frac{m}{2^s},\frac{n}{2^s}\Big)\in \mathcal I.
\end{equation}
\end{proposition}
$a\mid b$ means that $a$ divides $b$ for $a$ and $b$ integers. \\
{\bf Proof.}
We consider first $s=0$. For the case \eqref{4s} we have just $m^2+n^2$ is an odd number and that implies 
$(m,n)\in \mathcal J$.  Now assume $s>0$ so that $\la=4^s(2\ell+1)=m^2+n^2$.  
Since $\la$ is even we have 
either that both $m,n$ are odd or both $m,n$ are even. Suppose first that they are both odd. Then for some 
$(a,b)\in \mathbb N_0^2$ 
we have  $m=2a+1, \: n=2b+1$, $4^s(2\ell+1)=4(a^2+b^2+a+b)+2$ but this leads to 
$4^s(2\ell+1)/2=2(a^2+b^2+a+b)+1$, a contradiction. 

So we must have $m=2a,\:n=2b$. We get $4(a^2+b^2)=4^s(2\ell+1)$.  We can divide by 4 and obtain $a^2+b^2=4^{s-1}(2\ell+1)$. 
If $s=1$ then $(m/2,n/2)\in \mathcal J$. We proceed by recursion to obtain \eqref{J} in  proposition \ref{la=}.
\\
Now consider \eqref{2.4s}. We start with $s=0$. Then either  $m,n$ are both odd or both are even. If both are even, then 
$m=2a, \:n=2b$ and we get $4(a^2+b^2)=2(2\ell+1)$, a contradiction. Hence $(m,n)\in \mathcal I$. 
Now pick $s>0$. Then again both $m,n$ must be even and we can write  $m=2a,\: n=2b$ to get $a^2+b^2=2\cdot4^{s-1}(2\ell+1)$.    
We can proceed by recursion till we get the desired result. \finbox

Finally we use proposition \ref{la=} to show that any non-constant eigenfunction on the square cannot be stricly positive 
on the boundary.

We have 
\begin{equation}\label{mn}
u(x,0)=\sum_{(m,n)\in I(\la)}a_{m,n}(\cos mx+\cos nx),\:\: a_{m,n}=a_{n,m} 
\end{equation} 
and set without loss  $\sum_{(m,n)\in I(\la)}a_{m,n}=1$.\\\\
{\bf(1)} Assume first that $\la =2\ell+1$ so that  $(m,n)\in \mathcal J$. Then $\cos m\pi+\cos n\pi=0$ 
and hence $u(\pi,0)=0$.\\\\
{\bf(2)} If $(m,n)\in \mathcal I$ then
 $\cos m\pi+\cos n\pi=-2$.\\ We have  $u(\pi,0)=-2$ and $u(0,0)=u(2\pi,0)=2$.\\
{\bf(1a)} $\la=4^s(2\ell+1)$ with $s>1$. We have according to proposition \ref{la=} 
$\cos(2^{-s}m\pi)+\cos(2^{-s}n\pi)=0$ so that $u(2^{-s}\pi,0)=0$. \\  
{\bf(2a)} $\la=2\cdot 4^s(2\ell+1)$. Proposition \ref{la=} implies that $\cos(2^{-s}m\pi)+\cos(2^{-s}n\pi)=-2$.

This completes the proof of theorem \ref{I}. \finbox

\section{Some consequences}
If one takes theorem \ref{I} and its proof then it is easy to see that the following result are true.

\begin{corollary}\label{TorCy}

Take the torus $\mathbb T$, i.e. $\mathcal R(c,d)$ with periodic boundary conditions. Suppose  there is a real valued 
eigenfunction $u$  satisfying $-\Delta u=\la u$ with periodic boundary conditions, i.e. $u(0,y)=u(c,y)$, $u(x,0)=u(x,d)$,
$\pa_xu(0,y)=\pa_xu(c,y)$, $\pa_yu(x,0)=\pa_yu(x,d)$ with the property that $u>0$ on $\pa\mathcal R(c,d)$, 
then $u\equiv C$ for some positive constant $C$.\\ 
Similarly for the cylinder with Neumann boundary conditions which corresponds for $\mathcal R(c,d)$ with, say,
Neumann boundary conditions for $y=0$ and $y=d$ and for periodic boundary conditions for $x=0,\:x=c$ we have the same statement. 
If there is an eigenfunction $u$ satisfying those boundary conditions so that $u>0$ on $\pa \mathcal R(c,d)$ then
$u\equiv C$ for some positive constant $C$.
\end{corollary}
The proof is easy. Just note that if we have an eigenfunction that is positive on $\pa \mathbb R(c,d)$ and its expansion 
(for the torus in in the $x$ and the $y$ direction, for the cylinder just in the $x$ direction) contains also sines, we just can 
remove terms with the sines and keep the positivity and the boundary conditions. Hence we are back to the Neumann situation.    

Finally one should mention that for the isoceles triangle with a right angle,  say with with the vertices
$(0,0),(2\pi,0),(0,2\pi)$ we get the same result as for the square for the Neumann case since we can transform 
it to the problem for the square by symmetrization and reflection with respect to the line which goes through 
$(0,2\pi),(2\pi,0)$. 
\footnotesize
\bibliographystyle{plain}

\begin{thebibliography}{1}
\bibitem{CH1}
R. Courant, D. Hilbert.
\newblock Methods of Mathematical Physics, Volume 1, translated and revised from the German Original, 
Interscience Publishers, New York, 1953.
\bibitem{OW}
T. Hoffmann-Ostenhof.
\newblock Oberwolfach Reports. Report N0. 33/2012, Problem Section (xv).
European Mathematical Society Publishers, to appear. 
\bibitem{TothZelditch}
J. Toth, S. Zelditch
\newblock Counting nodal lines which touch the boundary of an analytic domain. 
\newblock arXiv:0710.0101
\end{thebibliography}

\scshape
\noindent
Thomas Hoffmann-Ostenhof: Department of Theoretical Chemistry,\\ A 1090 Wien, W\"ahringerstra\ss e 17, Austria\\
\noindent
email: thoffmann@tbi.univie.ac.at  
\end{document}